\documentclass{amsart}
\makeatletter
\@namedef{subjclassname@2020}{\textup{2020} Mathematics Subject Classification}
\makeatother
\usepackage{amsthm,amssymb,amsfonts,latexsym,mathtools,amsmath}
\usepackage[T1]{fontenc}
\usepackage{hyperref}
\hypersetup{
    colorlinks=true,
    linkcolor=blue,
    filecolor=blue,      
    urlcolor=blue,
    linktocpage=true
}
\newtheorem{theorem}{Theorem}[section]
\newtheorem{lemma}[theorem]{Lemma}
\newtheorem{proposition}[theorem]{Proposition}
\newtheorem{corollary}[theorem]{Corollary}
\theoremstyle{definition}

\numberwithin{equation}{section}
\newcommand{\kk}{\Bbbk}
\newcommand{\NN}{\mathbb{N}}
\newcommand{\ZZ}{\mathbb{Z}}
\newcommand{\Kalg}{\mathcal{K}}
\newcommand{\Oalg}{\mathcal{O}}
\DeclareMathOperator{\Oz}{Oz}
\DeclareMathOperator{\Aut}{Aut}
\DeclareMathOperator{\Frac}{Frac}
\DeclareMathOperator{\PIdeg}{PIdeg}
\DeclareMathOperator{\rank}{rank}
\DeclareMathOperator{\ord}{ord}
\DeclareMathOperator{\lcm}{lcm}
\DeclareMathOperator{\Spec}{Spec}

\newcommand{\paperhyphenation}{\lefthyphenmin=3\righthyphenmin=3\hyphenation{ozone discriminant arxiv graded ungraded groups}}

\newcommand{\paperabstract}{We determine centers of two families of four-dimensional double Ore extensions without torsion assumptions, restricting to $f\ne0$ in $\mathcal O$. On the PI locus of $\mathcal K$, we classify the ozone groups and obtain normal commutative fixed rings of generic degree two over the center, including a nonmonogenic example. For an odd-torsion subfamily of $\mathcal O$, independent odd orders $n\geq1,m>1$ yield ozone group $C_n^2\times C_2^2$, PI degree $2nm$ and fixed-algebra PI degree $m$, with unchanged center. A PI specialization of $\mathcal G$ has trivial ozone group.}

\begin{document}
\paperhyphenation
\title[Ozone groups and fixed rings of double Ore extensions]{Ozone groups and fixed rings of double Ore extensions}
\author{Andr\'es Rubiano}
\address{Universidad Distrital Francisco Jos\'e de Caldas}
\curraddr{Campus Universitario}
\email{aarubianos@udistrital.edu.co}
\date{}
\dedicatory{Dedicated to Abby}
\begin{abstract}
\paperabstract
\end{abstract}
\maketitle
\sloppy
\section{Introduction}

Let $\kk$ be an algebraically closed field of characteristic zero. For an associative unital $\kk$-algebra $R$, write $Z(R)=\{z\in R:zr=rz\text{ for every }r\in R\}$ for its center. If a group $H$ acts on $R$ by $\kk$-algebra automorphisms, its fixed algebra is $R^H=\{r\in R:h(r)=r\text{ for every }h\in H\}$.

The center of a PI algebra records its central characters, but it does not usually recover the algebra. Center-fixing automorphisms provide a second invariant of this relationship. For a $\kk$-algebra $A$, its ozone group is
\begin{equation}\label{eq:oz-definition}
 \Oz(A)=\{\phi\in\Aut_{\kk}(A):\phi|_{Z(A)}=\mathrm{id}\}.
\end{equation}
Chan, Gaddis, Won and Zhang introduced and developed this invariant in \cite{CGWZ24,CGWZ25}. In the noetherian PI Artin--Schelter regular setting, they proved finiteness and obtained central-rank bounds. Their characterization for algebras generated in degree one includes an abelian ozone-group hypothesis \cite[Theorem~F]{CGWZ25}. Liu, Wu and Zhu subsequently proved abelianness in this setting \cite{LWZ26}; together these results give $A^{\Oz(A)}=Z(A)$ if and only if $A$ is a skew polynomial ring. Determining the actual fixed algebra beyond that extremal case is the problem addressed here. Explicit calculations for further quantum algebras, including \cite{GomezVenegas26,GVB226,GSV26,WY26}, show the usefulness of identifying the entire center-fixing group rather than a selected subgroup of graded automorphisms.

Double Ore extensions provide a natural setting in which the relation between the center and the ozone fixed ring is more delicate. Zhang and Zhang introduced the construction in \cite{ZZ08} and used it to produce the regular families of type $(14641)$ in \cite{ZZ09}. The defining data may interchange generators instead of scaling each one separately. Consequently a diagonal calculation in the displayed generators need not exhaust the ozone group. Moreover, a power of an Ore automorphism can become inner in the quotient division ring before that automorphism becomes trivial. The corresponding central coefficients must be retained. The finite-order coefficient method in \cite[Lemma~2.3]{GKM17} and the extended-centroid analysis in \cite{Mat88} explain this distinction; our arbitrary-parameter center formulas are proved directly below.

The closest antecedent for the present center calculations is the author's earlier article, \emph{Center of double extension regular algebras of type $(14641)$} \cite{RubianoCenter26}. Its Theorem~3.7 distinguishes the complete centers in Table~1 from the central subalgebras in Table~2. For the normalized $\mathcal K$ presentation used there, Table~1 treats nontorsion parameters with $q=\pm1$, whereas Table~2 records central subalgebras at $f=\pm1$. Here we determine the complete center for arbitrary nonzero parameters, recover those nontorsion rows as a corollary and then compute the full center-fixing group and fixed ring on the three-parameter PI locus. The comparison in Section~\ref{sec:examples} explicitly identifies the mixed generators that complete those earlier central subalgebras.

For $\mathcal O$, the proof of \cite[Theorem~3.7(2)]{RubianoCenter26} establishes normality of $x_1^2-fx_2^2$ in the displayed $h_0=1$ presentation and obtains central powers. We retain $h_0$, determine the PI locus and describe the complete center for every allowed $f\ne0$ by signed exponent orbits. We then compute the ozone group and fixed algebra in the stated odd-order regime. Thus the earlier center computations supply concrete starting points for the present invariant-theoretic problem. Finite Gr\"obner--Shirshov bases and normal-ordering formulas provide further computational tools \cite{HHR25,RubianoNormal26}; explicit homological calculations for a $\mathcal K$ slice appear in \cite{RubianoHomological26}. Retaining the original scalars is essential: equivalence by a graded twist does not in general preserve the center as an algebra.

The preprints \cite{HHR25,RubianoCenter26,RubianoHomological26,RubianoNormal26} provide background and comparisons. All center descriptions, automorphism calculations and fixed-algebra arguments needed here are established below from the defining relations; no proof assumes a result from those four preprints.

The center formula treats the entire three-parameter family $\mathcal K$ without a torsion assumption; the group and fixed-ring calculations then cover its full PI locus. Its complete center is described by a central exponent lattice followed by an explicit action whose square is the identity. The ozone group is the union of an explicitly determined diagonal character group and one flip coset. If $M$ is the least common multiple of $2$ and the parameter orders, its order is $2M^2$, whereas the central rank is $4M^2$. Three $2$-adic cases determine the abstract group type. We identify its fixed ring as the center of the subalgebra of even degree in the two extension generators. This commutative fixed ring carries a nontrivial involution whose quotient is the original center. The resulting finite cover has generic degree two; it need not be globally generated by one quadratic element. We also recover the normalized parameters intrinsically from normal degree-one elements.

For $\mathcal O$, a general signed-orbit formula allows both common commutation signs, including $-1$ at even-order parameters. The ozone and fixed-algebra results concern a two-parameter odd-torsion subfamily. A simultaneous-swap Zhang twist gives a direct description of the center. Two independent odd orders $n\geq1$, $m>1$ determine
\begin{equation}\label{eq:intro-O}
 \Oz(A)\cong C_n^2\times C_2^2,\ 
 \PIdeg(A)=2nm,\ 
 \PIdeg(A^{\Oz(A)})=m.
\end{equation}
Here $Z(A^{\Oz(A)})=Z(A)$, so the fixed algebra is noncommutative. We also determine the PI locus of the full original $\mathcal O$ family; in particular its scalar slice $h_0=1$ contains no PI algebra. In a concrete specialization, the central spectrum is a product of two quadratic surface cones. The contribution is these uniform family calculations, including their residual fixed algebras. Low-order three-generator flip examples already occur in \cite[Section~8.2 and Question~6.7]{GKM17}. Quadratic invariant covers occur in \cite[Examples~5.4--5.5]{CGWZ25}; the Heisenberg calculation \cite[Propositions~4.9 and~4.12]{GSV26} recovers the center by a cyclic group acting on the ozone fixed ring. The formulas in \cite[Theorem~4.4 and Proposition~4.13]{GomezVenegas26} already allow a fixed abstract ozone group with varying PI degree. The specific contributions developed here are the complete $\mathcal K$ calculations, its nonmonogenic quadratic cover and the independent-order formulas for $\mathcal O$. The multigrading lemma, quantum-torus methods, finite-group invariant arguments and twisting construction retain their original attributions; reproducing their proofs does not change their provenance. Section~\ref{sec:scope} also exhibits a PI specialization of the original family $\mathcal G$ with trivial ozone group by adapting the $B_q(g)$ calculations of \cite{GaddisYee25,GomezVenegas26} to the original graded presentation.

Section~\ref{sec:preliminaries} records the ozone tools and the division-algebra interpretation of the residual fixed ring. Sections~\ref{sec:K} and~\ref{sec:O} treat $\mathcal K$ and $\mathcal O$, respectively. Section~\ref{sec:examples} gives explicit centers and compares them with \cite{RubianoCenter26}. Section~\ref{sec:scope} explains the choice of families and provides a specialization with trivial ozone group. Each formal statement specifies its algebra and parameter assumptions. Within each family section, $A$ always denotes the algebra explicitly defined in that statement; comparisons use $A_{\mathcal K}$ and $A_{\mathcal O}$. All algebra maps fix $\kk$ pointwise, $\NN$ includes zero, $C_r$ denotes a cyclic group of order $r$ and $\mu_r\subseteq\kk^\times$ denotes the $r$th roots of unity. The center lattices are denoted by $\Lambda_K$ and $\Lambda_O$; $\sigma_K$ is the Ore automorphism in the $\mathcal K$ family, whereas $\vartheta$ is the simultaneous swap in the $\mathcal O$ family. The even exponent for central powers in $\mathcal K$ is $N_K$, the intrinsic normal plane is $\mathcal V_K$ and the common odd exponent for $\mathcal O$ is $N_O$; the parameters of the final $\mathcal G$ specialization carry the subscript $G$.

\section{Ozone automorphisms and central quotients}\label{sec:preliminaries}

We use the double-extension convention of \cite[Definition~1.3]{ZZ08}, with the corrected formulation in \cite[Definition~1.1]{CLM11}. A trimmed double extension of $R$ is an algebra generated by $R,u,v$, free as a left $R$-module on $u^iv^j$, with
\begin{equation}\label{eq:DE}
 vu=p_{12}uv+p_{11}u^2,\ 
 \binom{u}{v}r=\sigma(r)\binom{u}{v},\quad r\in R,
\end{equation}
and with $p_{12}\ne0$, the corresponding right freeness and equality of the left and right spans of $u,v$. The matrices $\sigma(r)$ need not be triangular. The two main families are $\mathcal K$ and $\mathcal O$ of \cite[Section~4.3]{ZZ09}; Section~\ref{sec:scope} also treats a specialization of $\mathcal G$. Their noetherianity, domain property and Artin--Schelter regularity of global dimension four follow from \cite[Theorem~0.1]{ZZ09}. We also give direct PBW or twist realizations for the calculations.

For a nonzero normal element $b$ of a domain $A$ admitting a quotient division ring, our conjugation convention is
\begin{equation}\label{eq:eta}
 \eta_b(r)=brb^{-1},\  br=\eta_b(r)b.
\end{equation}
The inverse is taken in the quotient division ring. Normality ensures that $\eta_b$ restricts to an automorphism of the algebra and $\eta_b\in\Oz(A)$. This convention gives $\eta_{bc}=\eta_b\eta_c$.

The following is \cite[Lemma~1.9]{CGWZ25}, written for the multigrading used here and with our conjugation convention. We recall its proof to keep the applications explicit.

\begin{lemma}\label{lem:multigrading}
Let $\kk$ be an algebraically closed field of characteristic zero. Let $r\geq1$ be an integer and let $A=\bigoplus_{\nu\in\ZZ^r}A_\nu$ be a $\ZZ^r$-graded $\kk$-domain that is a finite module over its center $Z(A)$. Then every $\phi\in\Oz(A)$ satisfies $\phi(A_\nu)=A_\nu$ for every $\nu$. Moreover, there is a nonzero multihomogeneous normal element $b\in A$ such that $\phi(a)=bab^{-1}$ for all $a\in A$, with the inverse taken in the quotient division algebra of $A$.
In particular, $\Oz_{\mathrm{gr}}(A)=\Oz(A)$, where $\Oz_{\mathrm{gr}}(A)$ denotes the center-fixing automorphisms preserving every $A_\nu$.
\end{lemma}
\begin{proof}
Put $F=\Frac Z(A)$ and $\mathcal{D}=A\otimes_{Z(A)}F$.
This central localization is a finite-dimensional $F$-domain, hence a
division algebra. If $a/s\in Z(\mathcal{D})$, with
$a\in A$ and $0\ne s\in Z(A)$, then $a\in Z(A)$; consequently
$Z(\mathcal{D})=F$.
The automorphism $\phi$ extends to an $F$-automorphism of $\mathcal{D}$.
Skolem--Noether \cite[Tag~074R]{Stacks} gives $d\in\mathcal{D}^\times$ with
$\phi(a)=dad^{-1}$ for every $a\in A$.
Clearing a central denominator, write $d=b/s$ with $0\ne b\in A$ and
$0\ne s\in Z(A)$. Thus $ba=\phi(a)b$ for every $a\in A$.

Fix one coordinate grading and a nonzero multihomogeneous $a\in A$.
Since products of extreme homogeneous components do not vanish in a
domain, the largest degrees in this identity give
\begin{equation*}
 \max\deg b+\deg a=\max\deg\phi(a)+\max\deg b.
\end{equation*}
The same comparison of smallest degrees gives
$\min\deg\phi(a)=\deg a=\max\deg\phi(a)$.
Repeating for all coordinates proves
$\phi(A_\nu)\subseteq A_\nu$; applying the argument to $\phi^{-1}$ gives
equality. The identity $ba=\phi(a)b$ now separates into the
multihomogeneous components of $b$. Every nonzero component implements
$\phi$ and it is normal because $\phi$ is surjective.
\end{proof}

For a centrally finite domain $A$, put $F=\Frac Z(A)$ and $\mathcal{D}=A\otimes_{Z(A)}F$. As in the preceding proof, $\mathcal{D}$ is a finite-dimensional central division algebra. The definition of PI degree gives \cite[Chapter~13]{MR01}
\begin{equation}\label{eq:pidegree}
 \rank_{Z(A)}A=[\mathcal{D}:F]=\PIdeg(A)^2.
\end{equation}
The order/rank relation for ozone invariants is established in \cite[Proposition~3.9(1) and Theorem~E]{CGWZ25}. In the finite abelian setting, the next elementary lemma records its division-algebra form and factors the residual dimension through the center of the fixed division algebra.

\begin{lemma}\label{lem:residual}
Let $\kk$ be an algebraically closed field of characteristic zero. Let $F$ be a field extension of $\kk$, let $\mathcal{D}$ be a finite-dimensional central division algebra over $F$ and let $G$ be a finite abelian group acting faithfully on $\mathcal{D}$ by $F$-algebra automorphisms. Write $\mathcal{D}^G$ for its fixed division algebra. Then
\begin{equation}\label{eq:residual}
 [\mathcal{D}:F]=|G|\,[\mathcal{D}^G:F]
 =|G|\,\PIdeg(\mathcal{D}^G)^2[Z(\mathcal{D}^G):F].
\end{equation}
If, in addition, $A\subseteq \mathcal{D}$ is a $G$-stable $\kk$-domain with $F=\Frac Z(A)$ and $\mathcal{D}=A\otimes_{Z(A)}F$, then
$\mathcal{D}^G=A^G\otimes_{Z(A)}F$.
\end{lemma}
\begin{proof}
Decompose $\mathcal{D}$ into the character eigenspaces $\mathcal{D}_\chi$. Their support is a subgroup of the character group: products are nonzero and inverses exist. Faithfulness makes its annihilator trivial, so every character occurs. If $0\ne b\in \mathcal{D}_\chi$, then $\mathcal{D}_\chi=b\mathcal{D}^G$, since $b^{-1}c$ is invariant for every $c\in \mathcal{D}_\chi$. Hence $\mathcal{D}$ has dimension $|G|$ over $\mathcal{D}^G$. The second equality follows by taking the center of $\mathcal{D}^G$. Finally every element of $\mathcal{D}$ has a central denominator; an invariant fraction with such a denominator has invariant numerator.
\end{proof}

\section{The three-parameter family \texorpdfstring{$\mathcal K$}{Lg}}\label{sec:K}

For $a,q,h\in\kk^\times$, let $A=\Kalg(a,q,h)$ have generators $x,t,u,v$ and relations
\begin{equation}\label{eq:K-relations}
\begin{gathered}
 tx=qxt,\  vu=-uv,\  ux=axu,\  vx=axv,\\
 ut=tv,\  vt=htu.
\end{gathered}
\end{equation}
This is the full family $\mathcal K$ of \cite[Section~4.3, Algebra $\mathbb K$]{ZZ09}, in normalized coordinates. If the original common mixing scalar is $h_0$ and its remaining parameter is $f$, replace its second extension generator by $h_0$ times that generator; then $a=h_0$ and $h=h_0^2f$. In particular no scalar twist parameter has been discarded.

The presentation has two useful interpretations. Over $R=\kk_q[x,t]$, its double-extension matrix is
\begin{equation}\label{eq:K-DE}
 \sigma(x)=axI_2,\ 
 \sigma(t)=\begin{pmatrix}0&t\\ht&0\end{pmatrix},\  P=(-1,0).
\end{equation}
On the other hand, put
\begin{equation}\label{eq:K-Ore}
 \begin{gathered}
 B=\kk\langle x,u,v\rangle/(ux-axu,vx-axv,vu+uv),\\
 \sigma_K(x)=qx,\ \sigma_K(u)=h^{-1}v,\ \sigma_K(v)=u.
 \end{gathered}
\end{equation}
Then $A=B[t;\sigma_K]$, with $tb=\sigma_K(b)t$. This gives a PBW basis $x^iu^jv^kt^\ell$ and proves the domain property directly. Reordering $t$ gives the required left and right double-extension bases over $R$. Thus the off-diagonal presentation \eqref{eq:K-DE} does not mean that $A$ lies outside all iterated Ore presentations; compare \cite[Proposition~5.3]{ZZ09} and \cite{CLM11}.

\subsection{The PI locus and the complete center}

\begin{proposition}\label{prop:K-PI}
Let $\kk$ be an algebraically closed field of characteristic zero. Let $a,q,h\in\kk^\times$ and let $A=\Kalg(a,q,h)$ be the algebra defined by \eqref{eq:K-relations}. Then $A$ satisfies a polynomial identity if and only if $a,q,h$ are roots of unity. In this case $A$ is a finite module over its center $Z(A)$.
\end{proposition}
\begin{proof}
Suppose $A$ is PI. Its central quotient is finite dimensional by the central-quotient theorem for prime PI rings \cite[Chapter~13]{MR01}. In this quotient, a relation $rs=c sr$ between invertible elements gives $c^d=1$ by taking determinants of left multiplication, where $d$ is the central dimension. Apply this to the quantum-plane subalgebras on $(x,u)$, $(x,t)$ and $(u,t^2)$; their parameters are $a,q,h^{-1}$, respectively. Consequently each parameter is a root of unity.

Conversely choose an even integer $N_K$ divisible by the orders of $a,q,h$. The elements
\begin{equation*}
 x^{N_K},\  t^{2N_K},\  u^{N_K}+v^{N_K},\  u^{N_K}v^{N_K}
\end{equation*}
are central. The last two make $u^{N_K},v^{N_K}$ integral over their symmetric polynomial ring. The PBW basis therefore gives a finite spanning set over the central algebra generated by these four elements: first reduce the $x,t$ exponents, then use the quadratic equation for $u^{N_K},v^{N_K}$. Hence $A$ is centrally finite and PI.
\end{proof}

When $a,q,h$ are roots of unity, write
\begin{equation}\label{eq:M}
 M=\lcm\bigl(2,\ord(a),\ord(q),\ord(h)\bigr).
\end{equation}
The center construction below does not require torsion.

The subalgebra $C=\kk\langle x,u,v,w\rangle$, where $w=t^2$, is quantum affine space with the additional relations
\begin{equation}\label{eq:C-relations}
 wx=q^2xw,\  wu=h^{-1}uw,\  wv=h^{-1}vw.
\end{equation}
For $\nu=(i,j,k,\ell)\in\ZZ^4$, write $m_\nu=x^iu^jv^kw^\ell$ in its Laurent localization. Define the lattice
\begin{equation}\label{eq:L}
\Lambda_K=\left\{(i,j,k,\ell):
\begin{array}{l}
 a^{j+k}q^{2\ell}=1,\\
 a^i(-1)^kh^\ell=a^i(-1)^jh^\ell=1,\\
 q^{2i}h^{-(j+k)}=1
\end{array}\right\}.
\end{equation}

The commutation-character description of a quantum-affine center is
standard; in the PI case, see \cite[Lemma~4.1]{CGWZ24}.
For coefficient comparison in Ore-extension centers, compare
\cite[Lemma~2.3]{GKM17}. The latter is stated for a finite-order Ore
automorphism; the finite-support argument below handles arbitrary
nonzero $a,q,h$ directly.

\begin{theorem}\label{thm:K-center}
Let $\kk$ be an algebraically closed field of characteristic zero.
Let $a,q,h\in\kk^\times$ and let $A=\Kalg(a,q,h)$ be the algebra with
generators $x,t,u,v$ and relations \eqref{eq:K-relations}.
Set $w=t^2$ and $C=\kk\langle x,u,v,w\rangle\subseteq A$.
Let $\Lambda_K\subseteq\ZZ^4$ be the lattice in \eqref{eq:L} and write
$m_{(i,j,k,\ell)}=x^iu^jv^kw^\ell$ in the Laurent localization of $C$.
Let $\sigma_K$ be the automorphism of $C$ induced by conjugation by $t$,
so that $\sigma_K(x)=qx$, $\sigma_K(u)=h^{-1}v$,
$\sigma_K(v)=u$ and $\sigma_K(w)=w$.
Then
\begin{equation}\label{eq:K-center}
 Z(A)=Z(C)^{\langle\sigma_K\rangle},\ 
 Z(C)=\operatorname{span}_{\kk}\{m_\nu:\nu\in\Lambda_K\cap\NN^4\},
\end{equation}
where multiplication in $Z(C)$ is inherited from $C$ and
\begin{equation}\label{eq:tau-monomial}
 \sigma_K(m_{(i,j,k,\ell)})
 =q^ih^{-j}(-1)^{jk}m_{(i,k,j,\ell)}.
\end{equation}
The restriction of $\sigma_K$ to $Z(C)$ has square equal to the identity.
If, in addition, $a,q,h$ are roots of unity and
$M=\lcm(2,\ord(a),\ord(q),\ord(h))$, then
\begin{equation}\label{eq:K-rank}
 \PIdeg(A)=2M,\  \rank_{Z(A)}A=4M^2.
\end{equation}
\end{theorem}
\begin{proof}
Commuting a Laurent PBW monomial with $x,u,v,w$ gives exactly
\eqref{eq:L}. These commutation equations are diagonal on the monomial
basis, so they determine the entire center of $C$.
Conjugation by $t$ gives \eqref{eq:tau-monomial}; its square is
conjugation by $w$, hence is trivial on $Z(C)$.

The decomposition $A=C\oplus Ct$ separates the even and odd $t$-degrees.
If a nonzero odd term $ft$, with $f\in C$, were central, commuting it
with $u$ would give
\begin{equation*}
 uf=h^{-1}fv.
\end{equation*}
Since all commutation coefficients in $C$ are nonzero, the PBW supports
on the two sides are the support of $f$ translated by $e_u$ and $e_v$,
respectively. A finite nonempty subset of $\ZZ^4$ cannot be invariant
under translation by $e_u-e_v\ne0$.
Thus every central element belongs to $C$, where it must be central and
$\sigma_K$-invariant. The converse is immediate. This proves the center
formula for arbitrary nonzero $a,q,h$.

For the remaining claims, assume that $a,q,h$ are roots of unity.
Put $n_0=\ord(q^2)$ and
\begin{equation*}
 d_0=\lcm\bigl(2,\ord(a^{n_0}),\ord(h^{n_0})\bigr).
\end{equation*}
Use the unimodular Laurent coordinates $x,w,r=u^{-1}v,u$.
Here $r$ commutes with $x,w$ and anticommutes with $u$.
A monomial $x^iw^\ell r^ku^s$ is central exactly when
\begin{equation}\label{eq:index-conditions}
 (-1)^s=1,\quad q^{2\ell}a^s=1,\quad q^{-2i}h^s=1,\quad
 a^ih^\ell(-1)^k=1.
\end{equation}
The first three conditions force $s\in d_0\ZZ$.
At $s=0$, one has $i,\ell\in n_0\ZZ$ and the last condition has image
\begin{equation*}
 \langle a^{n_0},h^{n_0},-1\rangle=\mu_{d_0}.
\end{equation*}
Its kernel therefore has index $n_0^2d_0$ in the remaining three
coordinates. Conversely choose $i,\ell$ satisfying the middle
conditions at $s=d_0$. Then $(a^ih^\ell)^{d_0}=1$, since
$a^{d_0}=q^{-2\ell}$ and $h^{d_0}=q^{2i}$.
Changing $i,\ell$ by multiples of $n_0$ and changing $k$ removes the
last obstruction. The projection to the $s$-coordinate is consequently
exactly $d_0\ZZ$. Thus
\begin{equation}\label{eq:lattice-index}
 [\ZZ^4:\Lambda_K]=n_0^2d_0^2=M^2.
\end{equation}
For the last equality, use
$n_0\ord(c^{n_0})=\lcm(n_0,\ord(c))$ and
$2n_0=\lcm(2,\ord(q))$ to obtain
\begin{equation*}
 n_0d_0=\lcm\bigl(2n_0,\ord(a),\ord(h)\bigr)=M.
\end{equation*}

The Laurent quantum torus is free over its center on representatives
of the exponent cosets, so its division algebra $\mathcal D_K$ has
central dimension $M^2$. Write $F_K=Z(\mathcal D_K)$.
The division algebra of $A$ is $\mathcal D_K\oplus\mathcal D_Kt$:
this two-term crossed extension is a domain finite dimensional over
$\mathcal D_K$, hence is a division algebra.
The action of $t$ on $F_K$ is nontrivial, because it exchanges suitable
distinct powers of $u$ and $v$ and its square is inner.
An element of its center has no odd term, by commuting with $F_K$,
so its center is $F_K^{\sigma_K}$. Artin's fixed-field theorem
\cite[Tag~09I3]{Stacks} gives $[F_K:F_K^{\sigma_K}]=2$. Consequently
\begin{equation*}
 [\mathcal D_K\oplus\mathcal D_Kt:F_K^{\sigma_K}]
 =2M^2[F_K:F_K^{\sigma_K}]=4M^2.
\end{equation*}
Central localization identifies this field with $\Frac Z(A)$, proving
\eqref{eq:K-rank}.
\end{proof}

Formula \eqref{eq:K-center} also gives an explicit vector-space basis of the center. For $j>k$ take $m_\nu+\sigma_K(m_\nu)$, one for each two-element exponent orbit in $\Lambda_K\cap\NN^4$. For $j=k$ retain $m_\nu$ exactly when $q^ih^{-j}(-1)^j=1$. These orbit sums, with the inherited PBW multiplication, determine every central element, including mixed coefficients at smaller $t$-degree.

The following recovers the two nontorsion $\mathbb K$ entries in \cite[Theorem~3.7, Table~1]{RubianoCenter26}.

\begin{corollary}
\label{cor:K-nontorsion-centers}
Let $\kk$ be an algebraically closed field of characteristic zero,
let $h\in\kk^\times$ be not a root of unity, let $q\in\{1,-1\}$,
and let $A=\Kalg(1,q,h)$ be defined by \eqref{eq:K-relations}.
Then
\begin{equation*}
 Z(A)=
 \begin{cases}
  \kk[x],&q=1,\\
  \kk[x^2],&q=-1.
 \end{cases}
\end{equation*}
\end{corollary}
\begin{proof}
For a nonnegative vector $(i,j,k,\ell)\in\Lambda_K$, the last
condition in \eqref{eq:L} reads $h^{-(j+k)}=1$, since $q^2=1$.
As $h$ is not a root of unity, $j=k=0$.
The middle conditions then give $h^\ell=1$, hence $\ell=0$.
Thus $\Lambda_K\cap\NN^4=\{(i,0,0,0):i\geq0\}$ and $Z(C)=\kk[x]$.
Theorem~\ref{thm:K-center} identifies $Z(A)$ with its invariants under
$x\mapsto qx$, giving the two displayed rings.
\end{proof}

\subsection{The full ozone group}

For the remainder of this section assume that $a,q,h$ are roots of unity and retain $M$ from \eqref{eq:M}.

Let
\begin{equation}\label{eq:rho}
 \begin{gathered}
 \varpi_K:\ZZ^4\longrightarrow\ZZ^3,\ 
 \varpi_K(i,j,k,\ell)=(i,j+k,2\ell),\\
 \Gamma=\operatorname{Hom}(\ZZ^3/\varpi_K(\Lambda_K),\kk^\times).
 \end{gathered}
\end{equation}
Realize $\Gamma$ as the diagonal automorphisms
\begin{equation}\label{eq:Gamma-action}
 d_{\lambda,c,\beta}(x)=\lambda x,\quad
 d_{\lambda,c,\beta}(t)=\beta t,\quad
 d_{\lambda,c,\beta}(u,v)=c(u,v),
\end{equation}
with $\lambda^ic^{j+k}\beta^{2\ell}=1$ for all $(i,j,k,\ell)\in \Lambda_K$.

\begin{theorem}\label{thm:K-ozone}
Let $\kk$ be an algebraically closed field of characteristic zero. Let $a,q,h\in\kk^\times$ be roots of unity and let $A=\Kalg(a,q,h)$ be the algebra with generators $x,t,u,v$ and relations \eqref{eq:K-relations}. Put $M=\lcm(2,\ord(a),\ord(q),\ord(h))$, let $\Lambda_K$ be the lattice \eqref{eq:L} and define $\varpi_K(i,j,k,\ell)=(i,j+k,2\ell)$ and $\Gamma=\operatorname{Hom}(\ZZ^3/\varpi_K(\Lambda_K),\kk^\times)$, acting on $A$ by the diagonal maps \eqref{eq:Gamma-action}. Every $\kk$-algebra automorphism of $A$ fixing $Z(A)$ pointwise belongs to exactly one of the two cosets
\begin{equation}\label{eq:K-ozone}
 \Oz(A)=\Gamma\ \sqcup\ g\Gamma,\  g=\eta_t.
\end{equation}
The group is presented by $\Gamma$ and $g$, with
\begin{equation}\label{eq:K-extension}
 g\gamma=\gamma g\quad(\gamma\in\Gamma),\ 
 g^2=d_{q^2,h^{-1},1}\in\Gamma.
\end{equation}
In particular $|\Gamma|=M^2$ and $|\Oz(A)|=2M^2$. The conjugations induced by
\begin{equation}\label{eq:normal-generators}
 x,\quad t,\quad P_+=u^2+h^{-1}v^2,\quad
 P_-=u^2-h^{-1}v^2,\quad W=uv
\end{equation}
generate the entire group.
\end{theorem}
\begin{proof}
Give $A$ the free grading $\deg x=(1,0,0)$, $\deg t=(0,1,0)$ and $\deg u=\deg v=(0,0,1)$. Proposition~\ref{prop:K-PI} and Lemma~\ref{lem:multigrading} imply
\begin{equation*}
 \phi(x)=\lambda x,\ \phi(t)=\beta t,\ 
 \phi(u,v)=\mathsf M_\phi(u,v),\quad \mathsf M_\phi\in\mathrm{GL}_2(\kk).
\end{equation*}
Preserving $uv+vu=0$ forces $\mathsf M_\phi$ to be diagonal or antidiagonal: the coefficients of $u^2,v^2$ are twice the products of the two entries in their respective columns. The mixed relations then force $\mathsf M_\phi=cI$ or $\mathsf M_\phi=c\sigma_K|_{\kk u+\kk v}$. These formulas all define automorphisms.

Let $F_K=\Frac Z(C)$. A diagonal candidate restricts to an automorphism of $F_K$ fixing $F_K^{\sigma_K}=\Frac Z(A)$. By Artin's fixed-field theorem \cite[Tag~09I3]{Stacks}, its restriction is either the identity or $\sigma_K$. It cannot be $\sigma_K$, which exchanges distinct central powers of $u,v$, while a diagonal map scales each separately. It is the identity precisely when it annihilates $\Lambda_K$, which is \eqref{eq:Gamma-action}. Removing $g$ from an antidiagonal candidate gives the diagonal case. This proves exhaustion for all automorphisms, without assuming gradedness at the outset.

The kernel of $\varpi_K$ is generated by $(0,1,-1,0)$, whose intersection with $\Lambda_K$ is twice this subgroup. Also $\varpi_K(\ZZ^4)$ has index two in $\ZZ^3$. Equation \eqref{eq:lattice-index} gives
\begin{equation*}
 [\ZZ^3:\varpi_K(\Lambda_K)]=2\frac{[\ZZ^4:\Lambda_K]}{2}=M^2.
\end{equation*}
The square and commutation relations in \eqref{eq:K-extension} follow directly from the displayed actions and describe the group completely.

For the generator assertion, the actions of the five normal elements are
\begin{equation}\label{eq:normal-actions}
\begin{array}{c|ccc}
 b&\eta_b(x)&\eta_b(t)&\eta_b(u,v)\\ \hline
 x&x&q^{-1}t&a^{-1}(u,v)\\
 t&qx&t&(h^{-1}v,u)\\
 P_+&a^2x&ht&(u,v)\\
 P_-&a^2x&-ht&(u,v)\\
 W&a^2x&-ht&-(u,v).
\end{array}
\end{equation}
Their normality and these signs follow from \eqref{eq:K-relations}. Let $J$ be the diagonal subgroup generated by $\eta_x,\eta_{P_+}$, the pure signs obtained from $\eta_{P_-}\eta_{P_+}^{-1}$ and $\eta_W\eta_{P_-}^{-1}$ and $g^2$. Its annihilator on $\ZZ^3$, with coordinates $(i,s,d)$, is
\begin{equation*}
 a^sq^d=1,\quad a^{2i}h^d=1,\quad s,d\in2\ZZ,\quad q^{2i}h^{-s}=1.
\end{equation*}
These conditions are exactly $\varpi_K(\Lambda_K)$: write $d=2\ell$. The relation $a^{2i}h^{2\ell}=1$ makes $a^ih^\ell$ a sign. Since $s$ is even, choose $j+k=s$ with both parities prescribed by $(-1)^j=a^ih^\ell$. Finite character duality gives $J=\Gamma$. Adding $g$ proves the assertion.
\end{proof}

The lattice presentation records the exact action and its distinguished square; a Smith normal form of $\varpi_K(\Lambda_K)$ computes $\Gamma$ directly. The abstract isomorphism type has a simpler description in terms of parameter orders.

\begin{corollary}\label{cor:K-abstract}
Let $\kk$ be an algebraically closed field of characteristic zero. Let $a,q,h\in\kk^\times$ be roots of unity and let $A=\Kalg(a,q,h)$ be the algebra with generators $x,t,u,v$ and relations \eqref{eq:K-relations}. Set $M=\lcm(2,\ord(a),\ord(q),\ord(h))$. Let $v_2(d)$ denote the exponent of $2$ in a positive integer $d$ and put $v_a=v_2(\ord(a))$, $v_q=v_2(\ord(q))$ and $v_h=v_2(\ord(h))$. Then
\begin{equation}\label{eq:abstract-group}
 \Oz(A)\cong
 \begin{cases}
 C_{2M}\times C_M,&v_h\geq\max(v_a,v_q,1),\\
 C_M\times C_{M/2}\times C_2^2,&v_a>\max(v_q,v_h,1),\\
 C_M^2\times C_2,&\text{otherwise}.
 \end{cases}
\end{equation}
In particular relative primitive phases introduce no further abstract group type.
\end{corollary}
\begin{proof}
Write $G=\Oz(A)$, $g=\eta_t$, $\mathsf{X}=\eta_x$ and
$\mathsf{P}=\eta_{u^2+h^{-1}v^2}$. Let $\mathsf{S}$ be the sign on $t$ and let $\mathsf{R}$ be
the simultaneous sign on $u,v$. Theorem~\ref{thm:K-ozone} shows that
$G=\langle g,\mathsf{X},\mathsf{P},\mathsf{S},\mathsf{R}\rangle$ is abelian of order $2M^2$.
On the three blocks $x,t,(u,v)$, the diagonal actions of $\mathsf{X},\mathsf{P},g^2$ are
\begin{equation*}
\begin{array}{c|ccc}
 &x&t&(u,v)\\ \hline
 \mathsf{X}&1&q^{-1}&a^{-1}\\
 \mathsf{P}&a^2&h&1\\
 g^2&q^2&1&h^{-1}.
\end{array}
\end{equation*}
For a prime $p$, project $G$ to its $p$-primary subgroup and project
$a,q,h$ to their $p$-primary roots of unity. We use the same letters for
these projections. The displayed scalar actions remain valid.

First let $p$ be odd and put $d=p^{v_p(M)}>1$. The order of this primary
subgroup is $d^2$ and $\mathsf{S},\mathsf{R}$ project to the identity. Put $\mathsf{Q}=g^2$.
Each of $\mathsf{X},\mathsf{P},\mathsf{Q}$ has order dividing $d$.
At least one of $a,q,h$ has order $d$.
If $\ord(a)=d$, an equality $\mathsf{X}^i\mathsf{P}^j=1$ gives $a^{-i}=1$ from the
$u,v$ block and then $a^{2j}=1$ from the $x$ block. Since $d$ is odd,
these imply $d\mid i$ and $d\mid j$. Thus
$\langle \mathsf{X},\mathsf{P}\rangle\cong C_d^2$.
If $\ord(q)=d$, an equality $\mathsf{X}^i\mathsf{Q}^j=1$ gives $q^{-i}=1$ from $t$ and
$q^{2j}=1$ from $x$, again forcing $d\mid i,j$.
If $\ord(h)=d$, an equality $\mathsf{P}^i\mathsf{Q}^j=1$ gives $h^i=1$ from $t$ and
$h^{-j}=1$ from $u,v$, with the same conclusion.
In each case we have exhibited a subgroup of order $d^2$, so it is the
whole primary subgroup. The odd part of $G$ is therefore
$C_{M_{\rm odd}}^2$, where $M_{\rm odd}=M/2^{v_2(M)}$.

For the two-primary subgroup put $d=2^{v_2(M)}\geq2$; its order is
$2d^2$. The two-primary projection raises group elements to an odd power, so the projected $g$ still interchanges the lines $\kk u,\kk v$,
whereas $\mathsf{X},\mathsf{P},\mathsf{S},\mathsf{R}$ are diagonal. If $\ord(h)=d$, the $u,v$ scalar of
$g^2$ gives $\ord(g)=2d$, while the $t$ scalar gives $\ord(\mathsf{P})=d$.
If $g^i=\mathsf{P}^j$, comparison on $t$ gives $d\mid j$, hence
$g^i=\mathsf{P}^j=1$. Thus their cyclic subgroups intersect trivially and have
product of order $2d^2$. This is the first case of
\eqref{eq:abstract-group}.

Suppose next that $\ord(h)<d$ but $\ord(q)=d$.
Here $\ord(\mathsf{X})=d$. Moreover $g^2$ has order $d/2$, because its scalar
orders are $\ord(q^2)=d/2$ and $\ord(h)\mid d/2$; the nontrivial
flip parity of $g$ gives $\ord(g)=d$.
If $\mathsf{X}^ig^j=1$, the $t$ and $x$ coordinates give $d\mid i,j$,
respectively. The same two comparisons in an equality
$\mathsf{R}=\mathsf{X}^ig^j$ would force the right side to be the identity, a
contradiction. Hence
$\langle \mathsf{X},g,\mathsf{R}\rangle\cong C_d^2\times C_2$ and its order is the
whole $2d^2$.

It remains to consider $\ord(q),\ord(h)<d$.
If $\ord(a)=d$, write $q=a^{2b}$ and $h=a^{2c}$.
Then $\ord(\mathsf{X})=d$ and $\ord(\mathsf{P})=d/2$.
Indeed $\mathsf{X}^i\mathsf{P}^j=1$ forces $d\mid i$ by its $u,v$ scalar and
$d/2\mid j$ by its $x$ scalar. These comparisons also show that
$\mathsf{S}\notin\langle \mathsf{X},\mathsf{P}\rangle$. Consequently
\begin{equation*}
 \langle \mathsf{X},\mathsf{P},\mathsf{S}\rangle\cong C_d\times C_{d/2}\times C_2.
\end{equation*}
Define $\gamma_0=g\mathsf{X}^{-c}\mathsf{P}^{-b}$. Its square has scalar actions
\begin{equation*}
 q^2a^{-4b}=1\quad\hbox{on }x,\ 
 q^{2c}h^{-2b}=1\quad\hbox{on }t,\ 
 h^{-1}a^{2c}=1\quad\hbox{on }u,v.
\end{equation*}
Thus $\gamma_0^2=1$ and its nontrivial flip parity makes it independent
of the diagonal subgroup $\langle \mathsf{X},\mathsf{P},\mathsf{S}\rangle$. We obtain
$C_d\times C_{d/2}\times C_2^2$, of order $2d^2$.
For $d>2$, this is precisely the second case of
\eqref{eq:abstract-group}; for $d=2$, it is $C_2^3$ and belongs to the
third case. Finally, if all three projected parameters are trivial,
then $d=2$ and the independent involutions $\mathsf{S},\mathsf{R},g$ give $C_2^3$.
These alternatives exhaust the two-primary subgroup. Combining it
with the odd part gives exactly the three groups in
\eqref{eq:abstract-group}. The argument used only primitive orders and
the identities $q=a^{2b}$, $h=a^{2c}$ when needed, so relative primitive
phases yield no additional abstract type.
\end{proof}

\subsection{The parity center and its quadratic quotient}

Let $E$ be the subalgebra of even total degree in $u,v$. It is generated by $x,t,U=u^2,V=v^2,W=uv$, where $W^2=-UV$.

\begin{theorem}\label{thm:K-fixed}
Let $\kk$ be an algebraically closed field of characteristic zero. Let $a,q,h\in\kk^\times$ be roots of unity and let $A=\Kalg(a,q,h)$ be the algebra with generators $x,t,u,v$ and relations \eqref{eq:K-relations}. Let $E=\kk\langle x,t,u^2,v^2,uv\rangle\subseteq A$, equivalently the subalgebra of even total degree in $u,v$ and put $\mathcal F=A^{\Oz(A)}$. Then the exact ozone fixed ring is the commutative normal affine domain
\begin{equation}\label{eq:parity-center}
 \mathcal F=A^{\Oz(A)}=Z(E).
\end{equation}
It has a nontrivial involution $\kappa$ with
\begin{equation}\label{eq:double-cover}
 Z(A)=\mathcal F^{\langle\kappa\rangle},\ 
 [\Frac\mathcal F:\Frac Z(A)]=2.
\end{equation}
Thus $\Spec\mathcal F\to\Spec Z(A)$ is a finite quotient of normal affine varieties, of generic degree two.
\end{theorem}
\begin{proof}
Put $U=u^2$, $V=v^2$, $W=uv$ and $P_\pm=U\pm h^{-1}V$. The pure sign on $u,v$ belongs to the ozone group, so its invariants lie in $E$. The normal elements in \eqref{eq:normal-generators} generate $E$ as an algebra, since $P_+,P_-$ span $U,V$. An element of $E$ is fixed by their conjugations exactly when it commutes with those generators. Theorem~\ref{thm:K-ozone} therefore gives \eqref{eq:parity-center}.

Put $E_0=\kk\langle x,t^2,U,V,W\rangle$ and let $\sigma_K$ be conjugation by $t$ on $E_0$. The same finite-support comparison as in Theorem~\ref{thm:K-center}, now with $U,V$, excludes odd $t$-terms from $Z(E)$. Hence
\begin{equation*}
 \mathcal F=Z(E_0)^{\langle\sigma_K\rangle}.
\end{equation*}
Conjugation by $u$ restricts to the following polynomial automorphism of $E_0$:
\begin{equation}\label{eq:kappa-action}
 \theta(x)=ax,\quad \theta(t^2)=ht^2,\quad
 \theta(U)=U,\quad\theta(V)=V,\quad\theta(W)=-W.
\end{equation}
It commutes with $\sigma_K$ on these generators and $\theta^2=\eta_U$ is trivial on $Z(E_0)$. Thus it restricts to an involution $\kappa$ of $\mathcal F$. An element of $\mathcal F$ is central in $A$ exactly when it commutes with $u$: it already commutes with $x,t$ and $v=t^{-1}ut$ in the division algebra. This proves the ring equality in \eqref{eq:double-cover}.

The involution is nontrivial. Indeed
\begin{equation}\label{eq:Omega}
 \Omega=(U-h^{-1}V)W^{M-1}\in\mathcal F,\ 
 \kappa(\Omega)=-\Omega.
\end{equation}
It commutes with $x,t,U,V,W$ because $M$ is even and annihilates all parameter orders. It is nonzero by PBW independence. Multiplying a denominator by its conjugate identifies $(\Frac\mathcal F)^{\langle\kappa\rangle}$ with $\Frac Z(A)$. Artin's fixed-field theorem \cite[Tag~09I3]{Stacks} therefore gives the degree-two extension in \eqref{eq:double-cover}; finiteness of the quotient also follows from the module argument below. In particular $\Omega^2\in Z(A)$ and $\Frac\mathcal F=\Frac Z(A)(\Omega)$; this is a fraction-field statement.

For completeness, normality follows from a direct Laurent-center argument. Localizing $E_0$ at its monomials gives an even-parity quantum torus. Its center is a commutative Laurent polynomial ring on a sublattice of $\{(i,j,k,\ell)\in\ZZ^4:j+k\in2\ZZ\}$, in the original $(x,u,v,w)$ coordinates. Central powers identify the localization of $Z(E_0)$ with this Laurent center. The center of $E_0$ is the subring whose original exponent coordinates are nonnegative. If an element of its fraction field is integral over this subring, it first lies in the Laurent center; each coordinate minimum valuation is nonnegative, since a negative valuation would make the highest power the unique lowest term in a monic integral equation. Thus $Z(E_0)$ is normal. If a group $H$ acts on a normal domain $R$, an element of $\Frac(R^H)$ integral over $R^H$ is integral over $R$, hence lies in $R$ and, being invariant, in $R^H$. Applying this observation first to $Z(E_0)$ and then to $\mathcal F$ proves normality of $\mathcal F$ and $Z(A)$.

For affineness and finiteness, let $\mathcal P_K$ be the central polynomial algebra used in Proposition~\ref{prop:K-PI}. It lies in $Z(A)\subseteq Z(E_0)$ and $A$ is a finite module over the noetherian ring $\mathcal P_K$. Its $\mathcal P_K$-submodules $Z(E_0)$, $\mathcal F$ and $Z(A)$ are therefore finite as well, so these rings are affine. A finite set of $\mathcal P_K$-module generators of $\mathcal F$ also generates it over $Z(A)$, proving that the quotient in the statement is finite.
\end{proof}

\begin{corollary}\label{cor:K-isomorphism}
Let $\kk$ be an algebraically closed field of characteristic zero. Let $(a,q,h)$ and $(a',q',h')$ be triples of nonzero roots of unity in $\kk$ and let $A=\Kalg(a,q,h)$ and $A'=\Kalg(a',q',h')$ be defined by \eqref{eq:K-relations}, with their standard gradings in which all four generators have degree one. Then $A$ and $A'$ are isomorphic as graded $\kk$-algebras if and only if $(a,q,h)=(a',q',h')$.
\end{corollary}
\begin{proof}
The only normal degree-one lines are $\kk x$ and $\kk t$. To verify this, a normal degree-one element induces an ozone automorphism. Lemma~\ref{lem:multigrading} shows that each of its nonzero free-homogeneous components is normal and induces the same automorphism. A component in $\kk u+\kk v$ would be normal in the negative quantum plane, hence proportional to $u$ or $v$; neither is an eigenvector of $\sigma_K$, so neither is normal in $A$. The $x$ and $t$ components cannot both occur, because their conjugations on $\kk u+\kk v$ are, respectively, scalar and nonscalar.

Consequently the span $\mathcal V_K=\kk x+\kk t$ is intrinsic. On $A_1/\mathcal V_K$, conjugation by $x$ is the scalar $a^{-1}$ and conjugation by $t$ is nonscalar with square $h^{-1}I$. This distinguishes the two normal lines and recovers $a,h$. A graded isomorphism therefore sends $x,t$ to scalar multiples of $x',t'$ and their relation recovers $q$. The converse is immediate from the presentation.
\end{proof}

By Theorems~\ref{thm:K-center} and \ref{thm:K-ozone}, these algebras attain the largest proper divisor, one half, of the central rank. They are not skew polynomial rings, by \cite[Theorem~F]{CGWZ25}. Theorem~\ref{thm:K-fixed} identifies the residual algebra responsible for this factor: it is a quadratic field extension of the central fraction field.

\section{Centers and ozone invariants in \texorpdfstring{$\mathcal O$}{Lg}}\label{sec:O}

In the odd-order PI regime studied below, the ozone fixed ring is noncommutative. Let $A=\Oalg(h_0,f)$ be the original family $\mathcal O$ of \cite[Section~4.3.3]{ZZ09}, where $h_0\ne0$ and $f\ne1$:
\begin{equation}\label{eq:O-original}
\begin{gathered}
 x_2x_1=-x_1x_2,\  y_2y_1=-y_1y_2,\\
 y_1x_1=h_0(x_1y_1+fx_2y_2),\quad
 y_1x_2=h_0(-x_2y_1+x_1y_2),\\
 y_2x_1=h_0(fx_2y_1-x_1y_2),\quad
 y_2x_2=h_0(x_1y_1+x_2y_2).
\end{gathered}
\end{equation}
Here the exclusion is $f\ne1$, as in the original source; it also follows from the determinant $h_0^4(1-f)^2$ of its mixing matrix. We do not suppress $h_0$ by twist equivalence.

\subsection{The twist realization and its PI locus}

Suppose first that $f\ne0$, choose $\iota^2=-1$ and $e^2=f$ and put
\begin{equation}\label{eq:O-change}
 \begin{gathered}
 u=x_1-\iota e x_2,\quad v=x_1+\iota e x_2,\quad
 r=y_1+\iota y_2,\quad s=y_1-\iota y_2,\\
 \alpha=h_0(1+e),\ \beta=h_0(1-e).
 \end{gathered}
\end{equation}
Let $B$ be quantum affine space on $u,v,r,s$, with
\begin{equation}\label{eq:O-B}
 uv=vu,\quad rs=sr,\quad
 ru=\beta ur,\quad rv=\alpha vr,\quad
 su=\alpha us,\quad sv=\beta vs.
\end{equation}
The automorphism $\vartheta$ simultaneously exchanges $u,v$ and $r,s$. Its Zhang twist $B^\vartheta$ has the same graded vector space and multiplication
\begin{equation}\label{eq:twist-product}
 b*c=b\vartheta^{\deg b}(c)
\end{equation}
for homogeneous $b$; see \cite{Zhang96}. Direct substitution in \eqref{eq:O-original} gives
\begin{equation}\label{eq:O-binomial}
\begin{gathered}
 u*u=v*v,\quad r*r=s*s,\\
 r*u=\alpha v*s,\quad r*v=\beta u*s,\quad
 s*u=\beta v*r,\quad s*v=\alpha u*r.
\end{gathered}
\end{equation}
These are also the six defining relations obtained from \eqref{eq:O-B} by \eqref{eq:twist-product}. Thus the invertible change \eqref{eq:O-change} identifies $A=B^\vartheta$. The resulting graded surjection $A\to B^\vartheta$ is an isomorphism: both Hilbert series are $(1-z)^{-4}$, the first by double-extension freeness \cite{ZZ08,ZZ09} and the second by the definition of twisting \cite{Zhang96}. In particular the even subalgebras of $A$ and $B$ have identical multiplication. Finite permutation twists of skew polynomial rings also occur in the ozone computations of \cite[Section~2, Lemma~2.5 and Corollary~2.7]{CGWZ25}; the change of variables above gives the particular involution model used here for $\mathcal O$.

For the rest of this section juxtaposition in $u,v,r,s$ denotes multiplication in $B$; the symbol $*$ is retained when the $A$-product must be distinguished. Even elements and all the final central presentations can therefore be multiplied in $B$.

\begin{proposition}\label{prop:O-PI}
Let $\kk$ be an algebraically closed field of characteristic zero. Let $h_0\in\kk^\times$ and $f\in\kk\setminus\{1\}$ and let $A=\Oalg(h_0,f)$ be defined by \eqref{eq:O-original}. Choose $e\in\kk$ with $e^2=f$ and put $\alpha=h_0(1+e)$ and $\beta=h_0(1-e)$. Then $A$ is PI if and only if $f\ne0$ and $\alpha,\beta$ are roots of unity. On this locus $A$ is a finite module over its center $Z(A)$. In particular $\Oalg(h_0,0)$ is not PI for any $h_0\ne0$ and no member $\Oalg(1,f)$ with $f\ne1$ is PI.
\end{proposition}
\begin{proof}
If $f\ne0$ and $A$ is PI, its even subalgebra contains the quantum planes on $u^2,r^2$ and $v^2,r^2$, with parameters $\beta^4,\alpha^4$. The determinant argument in Proposition~\ref{prop:K-PI} makes both parameters torsion, hence makes $\alpha,\beta$ torsion.

Conversely, $B$ is finite over its monomial center when $\alpha,\beta$ are torsion. The even, $\vartheta$-invariant part of that center is central in $A$. For the finite group generated by parity and $\vartheta$, the orbit-polynomial argument makes $Z(B)$ integral over this invariant subring \cite[Tag~0BRE]{Stacks}. The finite algebra generating set of the affine ring $Z(B)$ also generates it over that subring, so the integral extension is a finite module extension. The invariant subring is affine by the Artin--Tate lemma \cite[Tag~00IS]{Stacks}. Multiplication by its elements is the same in $A$ and $B$, so $A$ is a finite module over this central subring and is PI \cite[Chapter~13]{MR01}.

At $f=0$, choose $\iota\in\kk$ with $\iota^2=-1$ and work directly with \eqref{eq:O-original}. Localize at the normal element $x_1$ and set $\Theta=x_2x_1^{-1}$, $r=y_1+\iota y_2$, $s=y_1-\iota y_2$. Calculation in $A$ gives
\begin{equation*}
 r\Theta=(-\Theta+\iota)r,\  s\Theta=(-\Theta-\iota)s,\ 
 (rs)\Theta=(\Theta-2\iota)(rs).
\end{equation*}
If $A$ were PI, these elements would lie in its finite-dimensional central division algebra of fractions \cite[Chapter~13]{MR01}. Conjugation would make left multiplication by $\Theta$ similar to left multiplication by $\Theta-2\iota$, contradicting equality of traces in characteristic zero.

Finally, when $h_0=1$, torsion $\alpha,\beta$ would satisfy $\alpha+\beta=2$. Embed their cyclotomic field into $\mathbb C$; equality in the triangle inequality forces $\alpha=\beta=1$, hence $f=0$, the excluded PI case.
\end{proof}

\subsection{The center for arbitrary parameters}

We compute the center of the displayed Zhang twist \cite{Zhang96} directly by PBW support and commutation characters. The signed-orbit formula below is the resulting calculation for \eqref{eq:O-original}.

\begin{theorem}\label{thm:O-general-center}
Let $\kk$ be an algebraically closed field of characteristic zero, let
$h_0\in\kk^\times$ and $f\in\kk\setminus\{0,1\}$ and let
$A=\Oalg(h_0,f)$ be the algebra with relations \eqref{eq:O-original}.
Choose $\iota,e\in\kk$ with $\iota^2=-1$ and $e^2=f$ and put
$\alpha=h_0(1+e)$ and $\beta=h_0(1-e)$.
Define $u,v,r,s$ by \eqref{eq:O-change}, let $B$ be the quantum affine
algebra \eqref{eq:O-B} and identify $A=B^\vartheta$, where
$\vartheta(u,v,r,s)=(v,u,s,r)$.
For $\nu=(i,j,k,\ell)\in\NN^4$ write
$b_\nu=u^iv^jr^ks^\ell$, $|\nu|=i+j+k+\ell$ and
$\vartheta\nu=(j,i,\ell,k)$ and define its commutation-character vector by
\begin{equation}\label{eq:O-character-vector}
 \boldsymbol\chi(\nu)=
 \bigl(\beta^k\alpha^\ell,\alpha^k\beta^\ell,
       \beta^{-i}\alpha^{-j},\alpha^{-i}\beta^{-j}\bigr).
\end{equation}
Thus $b_\nu x=\chi_x(\nu)xb_\nu$ for $x=u,v,r,s$.
For $\varepsilon\in\{1,-1\}$ set
\begin{equation}\label{eq:O-signed-exponents}
 \mathcal E_\varepsilon=
 \{\nu\in\NN^4:|\nu|\in2\NN,
       \ \boldsymbol\chi(\nu)=(\varepsilon,\varepsilon,\varepsilon,\varepsilon)\}.
\end{equation}
Then the complete center is
\begin{equation}\label{eq:O-general-center}
 Z(A)=\operatorname{span}_{\kk}
 \{b_\nu+\varepsilon b_{\vartheta\nu}:
          \varepsilon\in\{1,-1\},\ \nu\in\mathcal E_\varepsilon\}.
\end{equation}
Its multiplication is inherited from $B$.
More precisely, a basis consists of one displayed binomial from each
two-element $\vartheta$-orbit in each $\mathcal E_\varepsilon$, together
with $b_\nu$ for the fixed exponents $\vartheta\nu=\nu$ in $\mathcal E_1$.
Fixed exponents in $\mathcal E_{-1}$ contribute no central element.
No torsion assumption on $\alpha,\beta$ is required.
\end{theorem}
\begin{proof}
The formula for $\boldsymbol\chi(\nu)$ follows by commuting a PBW monomial
with the four generators of $B$. Since $uv=vu$ and $rs=sr$,
$\vartheta(b_\nu)=b_{\vartheta\nu}$ without an extra scalar.

There is no nonzero odd homogeneous central element. If such an element
had finite nonempty PBW support $\mathcal S$, its commutation with $u,v$
in the twisted multiplication would imply
\begin{equation*}
 \mathcal S+e_v=e_u+\vartheta(\mathcal S),\ 
 \mathcal S+e_u=e_v+\vartheta(\mathcal S).
\end{equation*}
Here every scalar introduced by PBW reordering is nonzero, so these are
equalities of supports. They force invariance of $\mathcal S$ under
the nonzero translation $2(e_u-e_v)$, which is impossible for a finite
nonempty subset of $\ZZ^4$.

For an even homogeneous element $z=\sum_\nu c_\nu b_\nu$, centrality is
equivalent to $zx=x\vartheta(z)$ for every $x=u,v,r,s$.
Comparing the coefficient of $xb_\nu$ gives
\begin{equation}\label{eq:O-center-coefficients}
 \chi_x(\nu)c_\nu=c_{\vartheta\nu}
 \ (x=u,v,r,s).
\end{equation}
Thus, whenever $c_\nu\ne0$, the four characters have a common value
$\varepsilon_\nu$. The action of $\vartheta$ permutes the first two and
last two characters, so the common value at $\vartheta\nu$ is also
$\varepsilon_\nu$. Applying \eqref{eq:O-center-coefficients} twice yields
$\varepsilon_\nu^2=1$. Consequently $\nu\in\mathcal E_{\varepsilon_\nu}$,
and $c_{\vartheta\nu}=\varepsilon_\nu c_\nu$.
For a two-element orbit this gives exactly one central binomial.
For a fixed orbit it gives $c_\nu=\varepsilon_\nu c_\nu$, which excludes
$\varepsilon_\nu=-1$ in characteristic zero.
Conversely, every listed binomial satisfies
\eqref{eq:O-center-coefficients}; it is therefore central.
The center is graded, so separating homogeneous components proves
\eqref{eq:O-general-center} and the basis assertion.
All its elements have even degree and hence their twisted product is
their ordinary product in $B$.
\end{proof}

For instance, take $\alpha=\iota$, $\beta=1$, so that
$h_0=(1+\iota)/2$ and $f=-1$.
The exponent $\nu=(6,2,2,2)$ belongs to $\mathcal E_{-1}$ and has a
two-element orbit. Hence
\begin{equation*}
 (u^6v^2-u^2v^6)r^2s^2\in Z(A).
\end{equation*}
This element anticommutes with every generator in $B$, so it does not
belong to $Z(B)$. In particular the negative-sign part cannot be discarded
outside the odd-order regime.

\subsection{Odd-order centers and ozone groups}

From now on suppose $\alpha,\beta$ have odd orders and
\begin{equation}\label{eq:nm}
 n=\ord(\alpha\beta)\geq1,\  m=\ord(\alpha/\beta)>1.
\end{equation}
There is no coprimality assumption on $n,m$. Write $\mathcal Z_O=Z(B)$ and $\mathcal Z_{O,\mathrm{ev}}$ for its even total-degree part. We use the center-lattice description for PI skew polynomial rings \cite[Lemma~4.1]{CGWZ24}. In the present parameters the complete center lattice of $B$ is the product of two copies of
\begin{equation}\label{eq:L0}
 \Lambda_O=\{(j_1,j_2)\in\ZZ^2:\beta^{j_1}\alpha^{j_2}=\alpha^{j_1}\beta^{j_2}=1\}.
\end{equation}
Thus $\mathcal Z_O$ is spanned by the PBW monomials whose two block exponents belong to $\Lambda_O$.

\begin{theorem}\label{thm:O-center}
Let $\kk$ be an algebraically closed field of characteristic zero. Let $h_0\in\kk^\times$, $f\in\kk\setminus\{0,1\}$ and let $A=\Oalg(h_0,f)$ be the algebra with generators $x_1,x_2,y_1,y_2$ and relations \eqref{eq:O-original}. Choose $e\in\kk$ with $e^2=f$ and set $\alpha=h_0(1+e)$, $\beta=h_0(1-e)$. Assume that $\alpha,\beta$ are roots of unity of odd order and that $\alpha\ne\beta$. Put $n=\ord(\alpha\beta)\geq1$ and $m=\ord(\alpha/\beta)>1$. Choose $\iota\in\kk$ with $\iota^2=-1$ and define $u,v,r,s$ by \eqref{eq:O-change}. Let $B$ be the quantum affine algebra \eqref{eq:O-B}, let $\vartheta(u,v,r,s)=(v,u,s,r)$ and set $\mathcal Z_O=Z(B)$ and $\mathcal Z_{O,\mathrm{ev}}=\bigoplus_{d\geq0}(\mathcal Z_O)_{2d}$. Identify $A$ with the Zhang twist $B^\vartheta$ through \eqref{eq:O-change}, so that these subspaces are viewed in its underlying graded vector space. Then
\begin{equation}\label{eq:O-center}
 Z(A)=\mathcal Z_{O,\mathrm{ev}}^{\langle\vartheta\rangle},\ 
 \rank_{Z(A)}A=4n^2m^2,\  \PIdeg(A)=2nm.
\end{equation}
\end{theorem}
\begin{proof}
Under the odd-order hypotheses every entry of \eqref{eq:O-character-vector} has odd order. Hence $\mathcal E_{-1}$ is empty and Theorem~\ref{thm:O-general-center} reduces to $Z(A)=Z(B)_{\mathrm{ev}}^{\langle\vartheta\rangle}$.

Put $N_O=\lcm(\ord(\alpha),\ord(\beta))$, which is odd. Modulo $N_O$, passage from $(j_1,j_2)$ to its sum and difference is invertible. Multiplying and dividing the two characters in \eqref{eq:L0} and using the invertibility of squaring in the odd root group, gives
\begin{equation*}
 (j_1,j_2)\in \Lambda_O\quad\Longleftrightarrow\quad
 j_1+j_2\equiv0\pmod n,\  j_1-j_2\equiv0\pmod m.
\end{equation*}
It follows that $[\ZZ^2:\Lambda_O]=nm$. In the associated quantum torus, PBW monomials representing the cosets of its center lattice form a basis over the center field, giving $\rank_{\mathcal Z_O} B=n^2m^2$. Parity and $\vartheta$ generate a faithful group $C_2^2$ on $\mathcal Z_O$: the independent central monomials $u^{N_O},v^{N_O},r^{N_O},s^{N_O}$ detect both actions. Multiplying a fraction's denominator by its other group conjugates identifies the fixed field with $\Frac Z(A)$. The fixed-field degree theorem for a faithful finite group action \cite[Tag~09I3]{Stacks} therefore gives $[\Frac \mathcal Z_O:\Frac Z(A)]=4$.

There is no rank-invariance assumption about arbitrary twists here. Every element of $Z(A)$ is even, $\vartheta$-invariant and $B$-central, so multiplication by it is literally identical in $A$ and $B$. These algebras are identical $Z(A)$-modules. Since $B$ is finite over $\mathcal Z_O$ and $\mathcal Z_O$ is finite over $Z(A)$, their common central rank is $4n^2m^2$, proving the remaining claims.
\end{proof}

\begin{theorem}\label{thm:O-ozone}
Let $\kk$ be an algebraically closed field of characteristic zero. Let $h_0\in\kk^\times$, $f\in\kk\setminus\{0,1\}$ and let $A=\Oalg(h_0,f)$ be the algebra with generators $x_1,x_2,y_1,y_2$ and relations \eqref{eq:O-original}. Choose $e\in\kk$ with $e^2=f$ and set $\alpha=h_0(1+e)$, $\beta=h_0(1-e)$. Assume that $\alpha,\beta$ are roots of unity of odd order and that $\alpha\ne\beta$. Put $n=\ord(\alpha\beta)\geq1$ and $m=\ord(\alpha/\beta)>1$. Then the full center-fixing $\kk$-automorphism group of $A$ is
\begin{equation}\label{eq:O-ozone}
 \Oz(A)\cong C_n^2\times C_2^2.
\end{equation}
Explicitly, its elements are indexed by $\xi,\eta\in\mu_n$, $\rho\in\{1,-1\}$ and $j\in\{0,1\}$ and act by
\begin{equation}\label{eq:O-action}
\begin{aligned}
 x_1&\longmapsto\rho\xi x_1,&
 x_2&\longmapsto(-1)^j\rho\xi x_2,\\
 y_1&\longmapsto\rho\eta y_1,&
 y_2&\longmapsto(-1)^j\rho\eta y_2.
\end{aligned}
\end{equation}
\end{theorem}
\begin{proof}
Use the variables $u,v,r,s$, the algebra $B$ and the swap $\vartheta$
from Theorem~\ref{thm:O-center}. Proposition~\ref{prop:O-PI} makes $A$
finite over its center. Lemma~\ref{lem:multigrading}, attributed to
\cite[Lemma~1.9]{CGWZ25} and applied to
$\deg x_i=(1,0)$ and $\deg y_i=(0,1)$, therefore shows that every
$\phi\in\Oz(A)$ preserves both degree-one blocks.

An invertible linear substitution preserving $x_2x_1=-x_1x_2$ is
diagonal or antidiagonal: comparison of the $x_1^2,x_2^2$ coefficients
forces the two column products to vanish. The same holds for $y_1,y_2$.
Write the two nonzero coefficients in the $x$-block as
$\lambda_1,\lambda_2$ and those in the $y$-block as
$\mu_1,\mu_2$, in that order. In a diagonal block these are the
coefficients of the original generators and in an antidiagonal block
they are the coefficients of the exchanged generators. Comparing PBW
coefficients in the first two mixed relations in
\eqref{eq:O-original} gives
\begin{equation*}
\begin{array}{c|c|l}
 x\text{-block}&y\text{-block}&\text{necessary equations}\\ \hline
 \text{diagonal}&\text{diagonal}&
 \lambda_1\mu_1=\lambda_2\mu_2,\quad
 \lambda_2\mu_1=\lambda_1\mu_2,\\
 \text{antidiagonal}&\text{antidiagonal}&
 \lambda_1\mu_1=f\lambda_2\mu_2,\quad
 \lambda_1\mu_2=f\lambda_2\mu_1,\\
 \text{diagonal}&\text{antidiagonal}&
 \lambda_1\mu_1=-\lambda_1\mu_1,\\
 \text{antidiagonal}&\text{diagonal}&
 \lambda_1\mu_1=-\lambda_1\mu_1.
\end{array}
\end{equation*}
The last two cases are impossible. In the first case
$\lambda_2=\delta\lambda_1$ and $\mu_2=\delta\mu_1$ for a common
$\delta\in\{1,-1\}$. In the second case
$\mu_2=\delta\mu_1$ and $\lambda_1=\delta f\lambda_2$.
The change of variables \eqref{eq:O-change} consequently gives precisely
\begin{equation*}
 \begin{gathered}
 \phi=D_{c,d,\epsilon}\vartheta^j,\ 
 D_{c,d,\epsilon}(u,v,r,s)=(cu,\epsilon cv,dr,\epsilon ds),\\
 c,d\in\kk^\times,\quad \epsilon\in\{1,-1\},\quad j\in\{0,1\}.
 \end{gathered}
\end{equation*}
Indeed, the first case gives $\epsilon=1$, $c=\lambda_1$,
$d=\mu_1$ and $(-1)^j=\delta$; the second gives $\epsilon=-1$,
$c=\iota\delta e\lambda_2$, $d=-\iota\mu_1$ and
$(-1)^j=-\delta$. Conversely these four forms preserve all six
relations \eqref{eq:O-binomial}. Under the present odd-order hypotheses,
$\vartheta$ fixes the center by \eqref{eq:O-center}; hence we may
remove $\vartheta^j$ and suppose $\phi=D_{c,d,\epsilon}$.

The two normal quadrics are
\begin{equation}\label{eq:O-quadrics}
 \mathfrak z_x=x_1*x_1-fx_2*x_2=uv,\ 
 \mathfrak z_y=y_1*y_1-y_2*y_2=rs.
\end{equation}
For $h_0=1$, the antecedent normal quadric $\mathfrak z_x$ and its
central powers appear in
\cite[Theorem~3.7(2), Table~2 and its proof]{RubianoCenter26};
those central powers generate the explicit central subalgebra recorded there.
Direct calculation from \eqref{eq:O-B} and \eqref{eq:twist-product}
gives
\begin{equation*}
 \begin{aligned}
 \eta_{\mathfrak z_x}(u,v,r,s)
   &=(u,v,(\alpha\beta)^{-1}r,(\alpha\beta)^{-1}s),\\
 \eta_{\mathfrak z_y}(u,v,r,s)
   &=((\alpha\beta)u,(\alpha\beta)v,r,s).
 \end{aligned}
\end{equation*}
Thus $\mathfrak z_x^n,\mathfrak z_y^n$ are central. Since
$D_{c,d,\epsilon}(\mathfrak z_x)=c^2\mathfrak z_x$ and
$D_{c,d,\epsilon}(\mathfrak z_y)=d^2\mathfrak z_y$, fixing these powers
forces
\begin{equation*}
 c^{2n}=d^{2n}=1.
\end{equation*}

Let $N_O$ be the least common multiple of the orders of $\alpha,\beta$;
it is odd and divisible by $n$. The nonzero element
$H_x=u^{2N_O}+v^{2N_O}$ is $B$-central, even and $\vartheta$-invariant,
so $H_x\in Z(A)$. If $\epsilon=-1$, then $D_{c,d,-1}$ multiplies
$u*v=u^2$ and $v*u=v^2$ by $-c^2$. These quadrics are even, so their
$A$-powers and $B$-powers coincide, yielding
\begin{equation}\label{eq:twist-sign}
 D_{c,d,-1}(H_x)=(-c^2)^{N_O}H_x=-c^{2N_O}H_x=-H_x.
\end{equation}
The last equality uses $c^{2n}=1$ and $n\mid N_O$.
This contradicts center-fixing in characteristic zero, so $\epsilon=1$.
This exclusion also applies when $n=1$.

Now $D_{c,d,1}$ commutes with $\vartheta$ and acts on $B$-monomials
by ordinary block scaling. The nonzero central element
\begin{equation*}
 u^{N_O}r^{N_O}+v^{N_O}s^{N_O}
\end{equation*}
is multiplied by $(cd)^{N_O}$. Thus $(cd)^{N_O}=1$.
Together with $(cd)^{2n}=1$ and $\gcd(2n,N_O)=n$, this implies
$(cd)^n=1$. Hence $c^n=d^n=\rho$ for a common sign and
uniquely $c=\rho\xi$, $d=\rho\eta$ with
$\xi,\eta\in\mu_n$, since $n$ is odd.

Conversely, the conjugations of the two quadrics generate the
independent $\mu_n$ block scalings. Global parity $\pi$ fixes the even
center, while $\vartheta$ fixes it by \eqref{eq:O-center}. These four
factors commute and are independent, since $n$ is odd. Translating
back through \eqref{eq:O-change} gives \eqref{eq:O-action} and proves
\eqref{eq:O-ozone} without an a priori finiteness or abelianness theorem.
For completeness, the normal elements $b_\pm=u^{N_O}\pm v^{N_O}$
satisfy $\eta_{b_+}=\vartheta$ and $\eta_{b_-}=\pi\vartheta$, directly
from \eqref{eq:twist-product}. Thus
$\mathfrak z_x,\mathfrak z_y,b_+,b_-$ implement generators of all four
factors.
\end{proof}

\subsection{The residual fixed algebra}

\begin{theorem}\label{thm:O-fixed}
Let $\kk$ be an algebraically closed field of characteristic zero. Let $h_0\in\kk^\times$, $f\in\kk\setminus\{0,1\}$ and let $A=\Oalg(h_0,f)$ be the algebra with generators $x_1,x_2,y_1,y_2$ and relations \eqref{eq:O-original}. Choose $e\in\kk$ with $e^2=f$ and set $\alpha=h_0(1+e)$, $\beta=h_0(1-e)$. Assume that $\alpha,\beta$ are roots of unity of odd order and that $\alpha\ne\beta$. Put $n=\ord(\alpha\beta)\geq1$ and $m=\ord(\alpha/\beta)>1$. Let $\mathcal F=A^{\Oz(A)}$. Then
\begin{equation}\label{eq:O-fixed}
 Z(\mathcal F)=Z(A),\  \PIdeg(\mathcal F)=m.
\end{equation}
In particular $\mathcal F$ is noncommutative and the two orders in \eqref{eq:nm} control the ozone group and the residual division algebra separately.
\end{theorem}
\begin{proof}
Use the realization $A=B^\vartheta$ from Theorem~\ref{thm:O-center} and put $N_O=\lcm(\ord(\alpha),\ord(\beta))$. Let $H$ act on $B$ by the two $\mu_n$ block scalings and global parity. Its invariant algebra $\mathcal F_0=B^H$ is spanned by the monomials in the lattice
\begin{equation}\label{eq:invariant-lattice}
 \Lambda_{\mathrm{inv}}=\{\nu\in\ZZ^4:S_x(\nu),S_y(\nu)\in n\ZZ,
 \ S_x(\nu)+S_y(\nu)\in2\ZZ\},
\end{equation}
restricted to nonnegative exponents. Here $S_x=i+j$, $D_x=i-j$, $S_y=k+\ell$, $D_y=k-\ell$. This algebra is even, so its multiplication is the same in $A$ and $B$ and Theorem~\ref{thm:O-ozone} gives $\mathcal F=\mathcal F_0^\vartheta$.

Choose an inverse $c$ of $2$ modulo the common odd order $N_O$. The cross-commutation character between the two blocks of $B$ is
\begin{equation*}
 (\alpha\beta)^{cS_xS_y}(\alpha/\beta)^{-cD_xD_y}.
\end{equation*}
On $\Lambda_{\mathrm{inv}}$ its first factor is one. The map $\Lambda_{\mathrm{inv}}\to(\ZZ/m\ZZ)^2$ given by $(D_x,D_y)$ is surjective: the vectors $(1,-1,0,0)$ and $(0,0,1,-1)$ belong to $\Lambda_{\mathrm{inv}}$ and map to $(2,0)$ and $(0,2)$; multiplication by two is invertible modulo $m$. Write $\omega=\alpha/\beta$. For $\nu,\mu\in \Lambda_{\mathrm{inv}}$, their monomials satisfy
\begin{equation*}
 b_\nu b_\mu=\omega^{c(D_x(\nu)D_y(\mu)-D_y(\nu)D_x(\mu))}b_\mu b_\nu.
\end{equation*}
This alternating pairing is nondegenerate on the two difference coordinates. Its radical therefore has index $m^2$ and consists precisely of $D_x,D_y\equiv0\pmod m$. Using the center-lattice method of \cite[Lemma~4.1]{CGWZ24}, these radical conditions, together with \eqref{eq:invariant-lattice}, are precisely the $B$-center conditions; coset monomials give the central rank. Hence
\begin{equation*}
 \begin{gathered}
 Z(\mathcal F_0)=\mathcal Z_{O,\mathrm{ev}},\  [\mathcal D_O:F_O]=m^2,
 \\
 F_O=\Frac \mathcal Z_{O,\mathrm{ev}},\  \mathcal D_O=\mathcal F_0\otimes_{\mathcal Z_{O,\mathrm{ev}}}F_O.
 \end{gathered}
\end{equation*}
The nonnegative lattice spans $\Lambda_{\mathrm{inv}}$, by adding sufficiently large even central powers, so the affine and Laurent calculations have the same rank.

The involution $\vartheta$ acts nontrivially on $F_O$, for example on $u^{2N_O},v^{2N_O}$. Its fixed field is $F_A=\Frac Z(A)$, by the same invariant-denominator argument used in Theorem~\ref{thm:O-center}. We use quadratic semilinear Galois descent \cite[Tag~0CDQ]{Stacks}, with its decomposition explicit here, to show that the fixed division algebra $\mathcal D_O^\vartheta$ has central dimension $m^2$ over $F_A$. Choose $\chi\in F_O$ with $\vartheta(\chi)=-\chi\ne0$. Every $b\in \mathcal D_O$ decomposes as
\begin{equation*}
 b=\frac{b+\vartheta(b)}2+
 \chi\frac{b-\vartheta(b)}{2\chi},
\end{equation*}
with both coefficients in $\mathcal D_O^\vartheta$. The sum $\mathcal D_O=\mathcal D_O^\vartheta\oplus \chi\mathcal D_O^\vartheta$ is direct over $F_A$, by applying $\vartheta$ and subtracting. Hence $[\mathcal D_O^\vartheta:F_A]=m^2$. An element central in $\mathcal D_O^\vartheta$ also commutes with the central element $\chi$ and therefore with this whole sum. It lies in $F_O\cap \mathcal D_O^\vartheta=F_A$, proving the asserted center.

Finally any invariant fraction with denominator $d\in \mathcal Z_{O,\mathrm{ev}}$ can be rewritten with the invariant denominator $d\vartheta(d)$. Thus $\mathcal D_O^\vartheta$ is precisely the central quotient of $\mathcal F$. Its center calculation gives $Z(\mathcal F)=\mathcal F\cap F_A=Z(A)$, because $A\cap\Frac Z(A)=Z(A)$. This proves \eqref{eq:O-fixed}.
\end{proof}

The integers $n,m$ can be prescribed independently as any odd integers $n\geq1$ and $m>1$. Choose odd-order roots $\lambda,\omega$ of those orders and take the unique square roots in their common odd root group with $\alpha^2=\lambda\omega$, $\beta^2=\lambda\omega^{-1}$ and $\alpha\beta=\lambda$. Then
\begin{equation*}
 h_0=\frac{\alpha+\beta}{2},\ 
 f=\left(\frac{\alpha-\beta}{\alpha+\beta}\right)^2
\end{equation*}
has $h_0\ne0$, $f\ne0,1$ and realizes \eqref{eq:nm}. Fixing $n$ while varying $m$ therefore keeps the group \eqref{eq:O-ozone} unchanged and realizes every odd residual PI degree $m>1$ within this double-extension family.

\section{Explicit centers and comparison}\label{sec:examples}

\subsection{The first two flip parameters}

Let $A_+=\Kalg(1,-1,1)$ and $A_-=\Kalg(1,-1,-1)$ over $\kk$, with the generators in \eqref{eq:K-relations}. In either algebra put $U=u^2,V=v^2,W=uv$.
For $h=1$, Theorem~\ref{thm:K-center} reduces to the involution
\begin{equation*}
 (x,U,V)\longmapsto(-x,V,U),\  t^2\longmapsto t^2.
\end{equation*}
Set $P=U+V$, $X=x^2$, $Y=x(U-V)$, $\Delta_K=(U-V)^2$ and $T_2=t^2$. Then
\begin{equation}\label{eq:K1-example}
 Z(A_+)
 \cong\kk[P,X,Y,\Delta_K,T_2]/(Y^2-X\Delta_K).
\end{equation}
Indeed the invariant monomials in $x,U-V$ have even total degree and their second Veronese ring has exactly these three generators and one relation. This is the standard simultaneous-sign invariant calculation (cf.\ \cite[Tag~0FK7]{Stacks}); Theorem~\ref{thm:K-center} shows that it gives the whole center. In particular the mixed element $x(U-V)$ must be retained alongside the familiar central powers and symmetric polynomials. In the notation above, the corresponding entry of \cite[Theorem~3.7, Table~2]{RubianoCenter26} is the central subalgebra $R_0=\kk[X,T_2,U+V,UV]$. Since $\Delta_K=(U+V)^2-4UV$, the complete center is $R_0[Y]/(Y^2-X((U+V)^2-4UV))$. The inclusion $R_0\subset Z(A_+)$ is strict: every element of $R_0$ has even $x$-degree, whereas the nonzero mixed element $Y=x(U-V)$ has odd $x$-degree.

Further invariant covers of centers already occur in \cite[Examples~5.4--5.5]{CGWZ25} and \cite[Propositions~4.9 and~4.12]{GSV26}. The next calculation establishes failure of generation by one element for the particular quadratic cover arising here.

\begin{proposition}\label{prop:K-nonmonogenic}
Let $\kk$ be an algebraically closed field of characteristic zero. Let $A=\Kalg(1,-1,1)$ be the algebra defined by \eqref{eq:K-relations}, let $Z=Z(A)$ be its center and let $\mathcal F=A^{\Oz(A)}$ be its full ozone fixed ring. Then $\mathcal F$ cannot be generated by one element as a $Z$-algebra.
\end{proposition}
\begin{proof}
Put $U=u^2$, $V=v^2$, $W=uv$, $P=U+V$, $Q=U-V$, $T_2=t^2$, $X=x^2$, $Y=xQ$ and $\Delta_K=Q^2$.
For $E_0=\kk\langle x,t^2,U,V,W\rangle$, the PBW basis gives
\[
 E_0=R\oplus RW,\  R=\kk[P,T_2,x,Q],\ 
 W^2=(Q^2-P^2)/4.
\]
Conjugation by $t$ fixes $P,T_2$ and negates $x,Q,W$.
Consequently $Z=Z(A)=\kk[P,T_2,X,Y,\Delta_K]/(Y^2-X\Delta_K)$ and, over
$C_0=\kk[P,T_2,X,\Delta_K]$, the ring $\mathcal F$ is free with basis
$1,Y,a_1=xW,a_2=QW$.
The relations $Ya_1=Xa_2$, $Ya_2=\Delta_Ka_1$ and
$W^2=(\Delta_K-P^2)/4$ show that at
$\mathfrak m=(P,T_2,X,Y,\Delta_K)\subset Z$ one has
\[
 \mathcal F/\mathfrak m\mathcal F
 \cong\kk[a_1,a_2]/(a_1^2,a_1a_2,a_2^2),
\]
with basis $1,a_1,a_2$. If $\mathcal F=Z[c]$, the usual orbit polynomial for the involution $\kappa$ from Theorem~\ref{thm:K-fixed} gives \cite[Tag~0BRE]{Stacks}
\[
 c^2-(c+\kappa(c))c+c\kappa(c)=0,
 \  c+\kappa(c),\ c\kappa(c)\in Z.
\]
Thus $1,c$ would generate $\mathcal F$ as a $Z$-module, forcing every
fiber to have dimension at most two. The displayed three-dimensional
fiber is a contradiction.
\end{proof}

For $A_-=\Kalg(1,-1,-1)$, the central subalgebra in \cite[Theorem~3.7, Table~2]{RubianoCenter26} uses $t^4$ instead of $t^2$. The additional mixed central element $Z_-=uvt^2$ reflects the inner-square phenomenon already present in the signed-swap prototype of \cite[Question~6.7]{GKM17}. With the same $P,X,Y,\Delta_K$, now put $T_4=t^4$. The complete center is
\begin{equation}\label{eq:Kminus-example}
 Z(A_-)\cong\frac{\kk[P,X,Y,\Delta_K,T_4,Z_-]}
 {(Y^2-X\Delta_K,\ Z_-^2+(P^2-\Delta_K)T_4/4)}.
\end{equation}
To verify completeness, the centralizer of $x,u,v$ is
\begin{equation*}
 \kk[x,U,V,T_4,Z_-]/(Z_-^2+UVT_4),
\end{equation*}
and conjugation by $t$ fixes $Z_-,T_4$ while sending $(x,U,V)$ to $(-x,V,U)$. Taking its invariants gives \eqref{eq:Kminus-example}. Both algebras have PI degree four and ozone order eight, but their ozone groups are $C_2^3$ and $C_4\times C_2$, respectively. These formulas extend the low-order Ore prototypes in \cite[Section~8.2 and Question~6.7]{GKM17} to the four-generator centers.

Even $M$ and the central rank together do not determine the abstract group throughout the full family. For $\iota^2=-1$, the actions \eqref{eq:normal-actions} give
\begin{equation*}
 \Oz\bigl(\Kalg(\iota,1,1)\bigr)\cong C_4\times C_2^3,
 \ 
 \Oz\bigl(\Kalg(1,1,\iota)\bigr)\cong C_8\times C_4.
\end{equation*}
Both have $M=4$. The first product follows from a common $u,v$-scaling of order four, independent $x$- and $t$-signs and the flip. The second follows from Corollary~\ref{cor:K-abstract}. This illustrates why Theorem~\ref{thm:K-ozone} specifies the character lattice and the distinguished square, rather than only the order.

\subsection{Two quadratic surface cones}

Let $\zeta\in\kk$ be primitive of odd order $\ell\geq3$. Take the algebra $A_\zeta=\Oalg(h_0,f)$ of \eqref{eq:O-original}, with
\begin{equation*}
 h_0=\frac{1+\zeta}{2},\ 
 f=\left(\frac{1-\zeta}{1+\zeta}\right)^2,
 \  n=m=\ell.
\end{equation*}
For the quantum affine algebra $B$ in \eqref{eq:O-B} with $\alpha=1$, $\beta=\zeta$, one has ${\mathcal Z_O}=Z(B)=\kk[{U_\ell},{V_\ell},{R_\ell},{S_\ell}]$, where ${U_\ell}=u^\ell,{V_\ell}=v^\ell,{R_\ell}=r^\ell,{S_\ell}=s^\ell$ are ordinary $B$-powers. Put
\begin{equation*}
 {P_+}={U_\ell}+{V_\ell},\quad {Q_+}={R_\ell}+{S_\ell},\quad {P_-}={U_\ell}-{V_\ell},\quad {Q_-}={R_\ell}-{S_\ell}.
\end{equation*}
Parity and $\vartheta$ become the independent signs on the pairs $({P_+},{Q_+})$ and $({P_-},{Q_-})$. Applying the standard second-Veronese invariant calculation to each pair (cf.\ \cite[Tag~0FK7]{Stacks}) gives
\begin{equation}\label{eq:O-cones}
 Z(A_\zeta)\cong
 \frac{\kk[p_0,p_1,p_2,q_0,q_1,q_2]}
 {(p_1^2-p_0p_2,\ q_1^2-q_0q_2)},
\end{equation}
with generators corresponding to ${P_+}^2,{P_+}{Q_+},{Q_+}^2,{P_-}^2,{P_-}{Q_-},{Q_-}^2$, all of degree $2\ell$. These expressions are computed in $B$ and viewed as even elements of $A_\zeta$; their mutual products agree in the two algebras. The center is thus a normal complete intersection whose spectrum is the product of two $A_1$ surface cones.

The central quotient map $\Spec {\mathcal Z_O}\to\Spec Z(A_\zeta)$ has degree four and group $C_2^2$. Its nonfree locus is exactly the union of the two planes ${P_+}={Q_+}=0$ and ${P_-}={Q_-}=0$. Their images are the two singular surfaces of \eqref{eq:O-cones}. Away from those planes the action is free, so the quotient is finite \emph{\'etale} \cite[Tag~0B87]{Stacks}. For one sign quotient the fiber at the cone vertex is $\kk[P,Q]/(P^2,PQ,Q^2)$, which is nonreduced; taking the product proves that both stated surfaces are branch loci. The Jacobian of the two cone equations identifies their union with the singular locus. These surfaces have codimension two, so the map is \emph{\'etale} in codimension one.

The fixed algebra is visibly noncommutative. In the $B$-notation, the two invariant elements
\begin{equation*}
 P_0=u^{2\ell-1}v+v^{2\ell-1}u,\ 
 Q_0=r^{2\ell-1}s+s^{2\ell-1}r
\end{equation*}
belong to $A_\zeta^{\Oz(A_\zeta)}$. Denote their two summands by $P_1,P_2$ and $Q_1,Q_2$. Their commutator is
\begin{equation*}
 [Q_0,P_0]=(\zeta^2-1)(P_1Q_1+P_2Q_2)
 +(\zeta^{-2}-1)(P_1Q_2+P_2Q_1)\ne0,
\end{equation*}
since the four PBW monomials are distinct and $\zeta^2\ne1$. Theorem~\ref{thm:O-fixed} gives the stronger exact value $\PIdeg(A_\zeta^{\Oz(A_\zeta)})=\ell$.

\section{Scope and further double Ore families}\label{sec:scope}

The families $\mathcal K$ and $\mathcal O$ were selected because their complete centers admit two explicit but different realizations. For $\mathcal K$, the decomposition $A=C\oplus Ct$ reduces the calculation to a quantum-affine central lattice and an involution. For $\mathcal O$, a simultaneous-swap Zhang twist gives the signed-orbit center formula, while odd-order commutation characters make the residual fixed division algebra explicit. These realizations account for the commutative quadratic cover in the first family and the independently prescribed residual PI degree in the second.

A uniform calculation for the remaining Zhang--Zhang labels would require separate parameter strata, PI criteria and proofs excluding every additional center-fixing automorphism. Complete-center descriptions or sufficiently strong central constraints must be established case by case. The central subalgebras in \cite[Theorem~3.7, Table~2]{RubianoCenter26} do not by themselves provide all those data: fixing a known central subalgebra is a weaker condition than fixing the full center. Actual PI skew-polynomial double extensions are already covered by \cite{CGWZ24,CGWZ25}. Other regimes require different arguments and the present results make no assertion that all uncomputed families are either difficult or trivial.

Small centers also behave differently from the PI cases studied here. If a nontrivial connected graded algebra generated in degree one has center $\kk$, every scalar homothety lies in its ozone group. Over our field this gives an infinite subgroup $\kk^\times$, rather than a trivial group. Thus an entry $Z(A)=\kk$ in a center table is not a trivial-ozone criterion.

\subsection{A specialization with trivial ozone group}

The next proposition gives a trivial ozone group inside the same 26-family list. It concerns one specified root-of-unity slice of $\mathcal G$; no conclusion about all parameters of $\mathcal G$ is intended. In this case a small collection of central elements already forces every center-fixing automorphism to be the identity, so an exhaustive affine-center presentation is unnecessary. The proof adapts the central-power and central-element calculations of \cite[Lemmas~3.1, 3.2, 3.8 and Proposition~3.9(3)]{GaddisYee25} and the ozone-automorphism argument of \cite[Lemma~4.1 and Proposition~4.3]{GomezVenegas26}. We verify the coordinate identification and the conclusion in the original graded algebra.

\begin{proposition}\label{prop:G-trivial}
Let $\kk$ be algebraically closed of characteristic zero, let ${p_G}\in\kk^\times$
have odd order $\ell\geq3$ and let $f\in\kk^\times$. Put ${q_G}={p_G}^{-1}$.
Let $A_G({p_G},f)$ be the graded $\kk$-algebra on degree-one generators
$c,X,Y,Z$ with $c$ central and relations
\begin{equation}\label{eq:G-trivial-relations}
 YZ={q_G}ZY,\  XY={q_G}YX-{q_G}^2cZ,\  XZ={q_G}^{-1}ZX-fcY.
\end{equation}
Then $A_G({p_G},f)$ is a PI member of the Zhang--Zhang family $\mathcal G$,
with original common mixing scalar $h_0={p_G}^{-1}$ and
$\Oz(A_G({p_G},f))=\{\mathrm{id}\}$.
\end{proposition}
\begin{proof}
In the presentation of \cite[Section~4.2.4, Algebra $\mathbb G$]{ZZ09}, set
$c=x_1$, $X=x_2+x_1/({p_G}-1)$, $Y=y_1$ and $Z=y_2$ and take $h_0={p_G}^{-1}$.
The relations become \eqref{eq:G-trivial-relations}.
They also present an Ore extension of
\begin{equation*}
 B_G=\kk[c]\otimes(\kk\langle Y,Z\rangle/(YZ-{q_G}ZY)),
\end{equation*}
with
$\sigma(Y)={q_G}Y$, $\sigma(Z)={q_G}^{-1}Z$,
$\delta(Y)=-{q_G}^2cZ$, $\delta(Z)=-fcY$ and $\delta(c)=0$.
Indeed $\delta(YZ-{q_G}ZY)=0$ and $\sigma(c)=c$.
Thus the monomials $c^iY^jZ^kX^r$ form a PBW basis and the algebra is a domain.

The following central-element formulas are the linear-polynomial adaptation of the cited $B_q(g)$ calculations \cite[Lemmas~3.1--3.2 and Proposition~3.9(3)]{GaddisYee25}. Put
\begin{equation*}
 a_0=\frac{f{q_G}}{{q_G}^2-1},\  b_0=\frac{{q_G}^4}{1-{q_G}^2},\ 
 \Omega_G=YZX+a_0cY^2+b_0cZ^2.
\end{equation*}
The relations show that $c,Y^\ell,Z^\ell,\Omega_G$ are central. For the powers,
use $[\ell]_{{q_G}^2}=[\ell]_{{q_G}^{-2}}=0$ in the identities
\begin{equation*}
 \begin{aligned}
 XY^r&={q_G}^rY^rX-{q_G}^{3-r}c[r]_{{q_G}^2}Y^{r-1}Z,\\
 XZ^r&={q_G}^{-r}Z^rX-fc[r]_{{q_G}^{-2}}YZ^{r-1}.
 \end{aligned}
\end{equation*}
Here $[r]_v=1+v+\cdots+v^{r-1}$. Commuting the displayed expression for
$\Omega_G$ with $X,Y,Z$ gives zero.
For the PI assertion we use the shifted quantum-torus method of \cite[Lemma~3.8 and Proposition~3.9(3)]{GaddisYee25}: take a central indeterminate $t$ and the quantum
division algebra over $\kk(c,t)$ generated by $Y,Z$ with $YZ={q_G}ZY$.
The assignment
\begin{equation*}
 X=(YZ)^{-1}(t-a_0cY^2-b_0cZ^2)
\end{equation*}
satisfies the defining relations and embeds $A_G({p_G},f)$ in that division
algebra: the highest $t$-degree of a PBW sum separates its largest
$X$-power and its remaining Laurent monomials are independent.
The target is finite dimensional over $\kk(c,t,Y^\ell,Z^\ell)$, so it is PI.

For the automorphism step, compare \cite[Lemma~4.1 and Proposition~4.3]{GomezVenegas26}. Write $A_G=A_G({p_G},f)$ and let $\phi\in\Oz(A_G)$. Fixing $Y^\ell,Z^\ell$
forces $\phi(Y),\phi(Z)$ to be homogeneous of degree one: compare highest
degrees in the domain, then constant terms. The quotient
$A_G/(Y,Z)=\kk[c,X]$ shows that their $c,X$ coefficients vanish.
Applying the same argument to $\phi^{-1}$ makes their restriction an
automorphism of the quantum plane on $Y,Z$. Comparing the $Y^2,Z^2,YZ$ coefficients in the preserved relation and using ${q_G}\ne\pm1$, forces the invertible linear map to be diagonal:
$\phi(Y)=\xi Y$, $\phi(Z)=\eta Z$, with $\xi^\ell=\eta^\ell=1$.
Also $\phi(c)=c$. Fixing $\Omega_G$ yields
\begin{equation*}
 YZ\bigl(\xi\eta\phi(X)-X\bigr)
 =a_0c(1-\xi^2)Y^2+b_0c(1-\eta^2)Z^2.
\end{equation*}
In PBW order $c,Y,Z,X$, every monomial on the left has positive exponents
in both $Y,Z$, whereas each monomial on the right has one exponent zero.
Hence $\xi^2=\eta^2=1$, because $a_0,b_0\ne0$. Oddness of $\ell$ gives
$\xi=\eta=1$ and cancellation then gives $\phi(X)=X$.
Thus $\phi=\mathrm{id}$, without any prior assumption that it is graded.
\end{proof}

This specialization is also related to the linear-polynomial calculations for $B_{q_G}(g)$ in \cite[Theorem~4.4]{GomezVenegas26}. After passing to the central field $\kk(c)$ and choosing $\rho_G^2=f$, the rescaling $Y'=\rho_G {q_G}^{-1}Y$, $Z'=Z$ gives $XY'={q_G}Y'X-{q_G}\rho_G cZ'$ and $XZ'={q_G}^{-1}Z'X-{q_G}\rho_G cY'$. The direct proof above establishes triviality for the original graded algebra without assuming that ozone groups transfer across a localization or a field extension.

\section{Conclusions and further questions}

The two families realize different residual algebras in Lemma~\ref{lem:residual}. For $\mathcal K$, the residual algebra is a quadratic field over the central fraction field. For the odd-torsion $\mathcal O$ family, it is a central division algebra of degree $m$ over its central fraction field. The factors left after taking ozone invariants are $2$ and $m^2$. For $A_{\mathcal K}=\Kalg(a,q,h)$ with the hypotheses of Theorem~\ref{thm:K-fixed} and $A_{\mathcal O}=\Oalg(h_0,f)$ with the hypotheses of Theorem~\ref{thm:O-fixed}, this reads
\begin{equation*}
 \frac{\rank_{Z(A_{\mathcal K})}A_{\mathcal K}}{|\Oz(A_{\mathcal K})|}=2,
 \ 
 \frac{\rank_{Z(A_{\mathcal O})}A_{\mathcal O}}{|\Oz(A_{\mathcal O})|}=m^2.
\end{equation*}
The complete center and the full ozone group enter both calculations. A central-power subalgebra alone would miss the inner-power coefficients in $\mathcal K$ and ordinary diagonal substitution would miss the sign obstruction in the Zhang twist. For $\mathcal K$, the abstract ozone group depends only on the parameter orders, whereas the graded algebra recovers the actual parameters.

Theorem~\ref{thm:O-general-center} also determines the center for even-order and nontorsion parameters, including the negative common character. Their ozone groups and residual fixed algebras are not determined by the odd-order calculation and form a concrete next problem. More broadly, it is natural to identify which other double-extension families admit a parity subalgebra whose center recovers the ozone fixed ring and when the residual fixed algebra is a field rather than a central division algebra of higher degree.

\section*{Acknowledgments}

The author is grateful to Professor Claudia Gallego and to the members of the Algebra Seminar, in particular Helbert Venegas, H\'ector Su\'arez and James G\'omez, for their valuable contributions, insightful comments and helpful mathematical discussions.

\end{document}